\documentclass[12pt]{article}
\usepackage{latexsym,amssymb,amsmath}
\begin{document}
\baselineskip=20pt

\newcommand{\la}{\langle}
\newcommand{\ra}{\rangle}
\newcommand{\psp}{\vspace{0.4cm}}
\newcommand{\pse}{\vspace{0.2cm}}
\newcommand{\ptl}{\partial}
\newcommand{\dlt}{\delta}
\newcommand{\sgm}{\sigma}
\newcommand{\al}{\alpha}
\newcommand{\be}{\beta}
\newcommand{\G}{\Gamma}
\newcommand{\gm}{\gamma}
\newcommand{\vs}{\varsigma}
\newcommand{\Lmd}{\Lambda}
\newcommand{\lmd}{\lambda}
\newcommand{\td}{\tilde}
\newcommand{\vf}{\varphi}
\newcommand{\yt}{X^{\nu}}
\newcommand{\wt}{\mbox{wt}\:}
\newcommand{\rd}{\mbox{Res}}
\newcommand{\ad}{\mbox{ad}}
\newcommand{\stl}{\stackrel}
\newcommand{\ol}{\overline}
\newcommand{\ul}{\underline}
\newcommand{\es}{\epsilon}
\newcommand{\dmd}{\diamond}
\newcommand{\clt}{\clubsuit}
\newcommand{\vt}{\vartheta}
\newcommand{\ves}{\varepsilon}
\newcommand{\dg}{\dagger}
\newcommand{\tr}{\mbox{Tr}}
\newcommand{\ga}{{\cal G}({\cal A})}
\newcommand{\hga}{\hat{\cal G}({\cal A})}
\newcommand{\Edo}{\mbox{End}\:}
\newcommand{\for}{\mbox{for}}
\newcommand{\kn}{\mbox{ker}}
\newcommand{\Dlt}{\Delta}
\newcommand{\rad}{\mbox{Rad}}
\newcommand{\rta}{\rightarrow}
\newcommand{\mbb}{\mathbb}
\newcommand{\lra}{\Longrightarrow}

\begin{center}{\Large \bf Differential Equations for Singular Vectors of $sl(n)$ }\end{center}
\pse

\begin{center}{\large Xiaoping Xu}\end{center}
\begin{center}{Institute of Mathematics, Academy of Mathematics \& System Sciences}\end{center}
\begin{center}{Chinese Academy of Sciences, Beijing 100080, P.R.China}\end{center}

\vspace{0.3cm}

\begin{center}{\Large \bf Abstract}\end{center}
\vspace{0.2cm}

{\small Given a weight of $sl(n)$, we derive a system of variable-coefficient second-order linear partial differential equations that determines the singular vectors in the corresponding Verma module. Moreover, we completely solve the system in a certain space of power series. The polynomial solutions correspond to the singular vectors in the Verma module. The well-known results of Verma, Bernstein-Gel'fand-Gel'fand and Jantzen for the case of $sl(n)$ are naturally included in our almost elementary approach of partial differential equations.} 

\section{Introduction}

One of the most beautiful things in Lie algebras is the highest weight representation theory. It was established based on the induced modules of  a Lie algebra with respect to a Cartan decomposition  from one-dimensional modules  of the  Borel subalgebra associated with a linear function (weight) on the Cartan subalgebra. These modules are now known as {\it Verma modules} [V1]. A {\it singular vector} (or {\it canonical vector}) in a Verma module is an invariant vector under the action of the Borel subalgebra. It is well known that the structure of a Verma module of a finite-dimensional simple Lie algebra is completely determined by its singular vectors (cf. [V1]). In this paper, we find  explicit formulas for the singular vectors in  the Verma modules of the Lie algebra $sl(n)$, in terms of differential operators.

 The structure of Verma module was first studied by Verma [V1].  Verma reduced the problem of determining all submodules of a Verma module of a finite-dimensional  semisimple Lie algebra to determining the embeddings of the other Verma modules into the objective module. He  proved that the multiplicity of the embedding is at most one. Bernstein, Gel'fand and Gel'fand [BGG] introduced the well-known useful notion of category ${\cal O}$ of representations, and found a necessary and sufficient condition for the existence of such an embedding in terms of the action of  Weyl group on weights. Sapovolov [S] introduced a certain bilinear form on a universal enveloping algebra. Lepowsky [L1-L4] studied analogous induced modules with respect to Iwasawa decomposition that is more general than Cartan decomposition, and obtained similar results as those in [V1] and [BGG]. These modules are now known as {\it generalized Verma modules}.

 Jantzen [J1,J2] introduced his famous "Jantzen filtrations" on Verma modules and used Sapovolov form to determine weights of singular vectors in Verma modules.  Verma modules of infinite-dimensional Lie algebras were first studied by Kac [K1]. Kac and Kazhdan [KK] generalized the results of Verma [V1] and Bernstein-Gel'fand-Gel'fand [BGG]  to the contragredient Lie algebra corresponding a symmetrizable generalized Cartan matrix. Deodhar, Gabber and Kac [DGK] generalized the results further to more general matrix Lie algebras. Rocha-Caridi and Wallach [RW1, RW2] generalized the results of Verma [V1] and Bernstein-Gel'fand-Gel'fand [BGG] to a class of graded Lie algebra possessing a Cartan decomposition and obtained Jantzen's character formula corresponding to the quotient of two Verma modules. The resolutions of irreducible highest weight modules over rank-2 Kac-Moody algebras were constructed. 

     One of the fundamental and difficult remaining problems in this direction is how to determine the singular vectors explicitly. Malikov, Feigin and Fuchs [MFF] introduced a formal manipulation on products of several general powers of negative simple root vectors and used free Lie algebras to give a rough condition when such a product is well defined. It seems to us that their condition can not be verified in general and their method can practically be applied only to finding very special singular vectors.

     In this paper, we introduce an almost elementary partial differential equation approach of determining the singular vectors in any Verma module of $sl(n)$. First, we identify the Verma modules with a space of polynomials, and the action of $sl(n)$ on the Verma module is identified with a differential operator action of $sl(n)$ on the polynomials. Any singular vector in the Verma module becomes a polynomial solution of a system of variable-coefficient second-order linear partial differential equations. Thus we have changed a difficult problem in a noncommutative space to a problem in commutative space. However, it is impossible to solve the system in the space of polynomials. So we extend the action of $sl(n)$ on the polynomial space to a larger space of certain formal power series. On this larger space, the negative simple root vectors become differential
operators whose arbitrary complex powers are well defined (so are their products). In this way, we overcome the difficulty of determining whether a product of several general powers of negative simple root vectors is well defined in the work [MFF] of Malikov, Feigin and Fuchs. 

Using commutator relations among root vectors and a certain substitution-of-variable technique that we developed in [X], we completely solve the system of partial differential equations in the space of power series. There are exactly $n!$ linearly independent solutions, where the polynomial solutions correspond to singular vectors in the Verma module. The Verma module has $n!$ linearly independent singular vectors if and only if the weight is dominant integral. 

As an application, we give three examples on the structure of the Verma modules of $sl(n)$ with certain special highest weights that are not dominant integral. In a subsequent work, we will show that the first two examples  give rise to generalizations of the wedge representations of the Lie algebra $W_{1+\infty}$ (cf. {KP], [KR1]). Moreover, the third example  will imply that the vacuum representations of  $W_{1+\infty}$ with negative integral levels are indeed irreducible.

 The result in this paper is not just useful in determining the structure of infinite-dimensional highest weight irreducible representations, but also give a new way of using symmetry to solve systems of partial differential equations.

Throughout this paper, all the vector spaces (algebras) are assumed over $\mbb{C}$, the field of complex numbers. Denote by $\mbb{Z}$ the ring of integers and by $\mbb{N}$ the additive semigroup of nonnegative integers.

In Section 2, we will derive the system of partial differential equations and give certain exact solutions. In Section 3, we will completely solve the system in a certain space of power series. In Section 4, we give three examples of application.

\section{Equations and Exact Solutions}

In this section, we will first derive a system of variable-coefficient second-order partial differential equations that determines the singular vectors in the Verma modules over the special Lie algebra  ${\rm sl}(n)$. Then we will present and prove certain exact solution of the system.

The special  Lie algebra
$$sl(n)=\mbox{Span}\{E_{i,j},E_{j,i}, E_{i,i}-E_{i+1,i+1}\mid 1\leq i<j\leq n\},\eqno(2.1)$$
with the Lie bracket:
$$[A,B]=AB-BA\qquad \for\;\;A,B\in sl(n),\eqno(2.2)$$
where $E_{i,j}$ is the $n\times n$ matrix whose $(i,j)$-entry is 1 and the others are 0. Set

$$h_i=E_{i,i}-E_{i+1,i+1},\qquad i=1,2,...,n-1.\eqno(2.3)$$
The subspace
$$H=\sum_{i=1}^{n-1}\mbb{C}h_i\eqno(2.4)$$
forms a Cartan subalgebra of $sl(n)$. We choose
$$\{E_{i,j}\mid 1\leq i<j\leq n\}\;\;\mbox{as positive root vectors}.\eqno(2.5)$$
In particular, we have 
$$\{E_{i,i+1}\mid i=1,2,...,n-1\}\;\;\mbox{as  positive simple root vectors}.\eqno(2.6)$$
Accordingly, 
$$\{E_{i,j}\mid 1\leq j<i\leq n\}\;\;\mbox{are negative root vectors}\eqno(2.7)$$
and
$$\{E_{i+1,i}\mid i=1,2,...,n-1\}\;\;\mbox{are negative simple root vectors}.\eqno(2.8)$$

Let 
$$\G=\sum_{1\leq j<i\leq n}\mbb{N}\es_{ij}\eqno(2.9)$$
be the torsion-free additive semigroup  of  rank $n(n-1)/2$ with $\es_{ij}$ as base elements, and let $U(sl(n))$ be the universal enveloping algebra of $sl(n)$. For 
$$\al =\sum_{1\leq j<i\leq n}\al_{ij}\es_{ij}\in \G,\eqno(2.10)$$
we denote
$$E^\al=E_{2,1}^{\al_{2,1}}E_{3,1}^{\al_{3,1}}E_{3,2}^{\al_{3,2}}E_{4,1}^{\al_{4,1}}\cdots E_{n,1}^{\al_{n,1}}\cdots E_{n,n-1}^{\al_{n,n-1}}\in U(sl(n)).\eqno(2.11)$$
Denote by  ${\cal G}_-$  the Lie subalgebra spanned by (2.7)  and by $U({\cal G}_-)$ its  universal enveloping algebra. Then 
$$\{E^\al\mid \al\in\G\}\;\;\mbox{forms a PBW basis of}\;\;U({\cal G}_-).\eqno(2.12)$$

 Let $\lmd$ be a weight, which is a linear function
on $H$, such that
$$\lmd(h_i)=\lmd_i\qquad\for\;\;i=1,2,...,n-1.\eqno(2.13)$$
Recall that $sl(n)$ is generated by $\{E_{i,i+1},E_{i+1,i}\mid i=1,2,...,n-1\}$ as a Lie algebra. The Verma  $sl(n)$-module with the highest-weight vector $v_\lmd$ of weight $\lmd$ is given by
$$M_\lmd =\mbox{Span}\{E^\al v_\lmd\mid\al\in \G\},\eqno(2.14)$$
with the action determined by
\begin{eqnarray*}E_{i,i+1}(E^\al v_\lmd)&=&(\sum_{j=1}^{i-1}\al_{i+1,j}E^{\al+\es_{i,j}
-\es_{i+1,j}}-\sum_{j=i+2}^n\al_{j,i}E^{\al+\es_{j,i+1}-\es_{j,i}}\\ & &+\al_{i+1,i}(\lmd_i+1-\sum_{j=i+1}^n\al_{j,i}+\sum_{j=i+2}^n\al_{j,i+1})E^{\al-\es_{i+1,i}})v_\lmd,\hspace{2.3cm}(2.15)\end{eqnarray*}
$$E_{i+1,i}(E^\al v_\lmd)=(E^{\al+\es_{i+1,i}}+\sum_{j=1}^{i-1}\al_{i,j}E^{\al+\es_{i+1,j}-\es_{i,j}})v_\lmd\eqno(2.16)$$
for $i=1,...,n-1$. For any $\al\in\G$, we define the {\it weight} of $E^\al v_\lmd $  by
$$({\rm wt}\:E^\al v_\lmd )(h_i)=(\lmd_i+\sum_{p=1}^{i-1}(\al_{i,p}-\al_{i+1,p})+\sum_{j=i+2}^n(\al_{j,i+1}-\al_{j,i})-2\al_{i+1,i})h_i\eqno(2.17)$$
 for $i=1,...,n-1$. Then the Verma module $M_\lmd$ is a space graded by weights.  A {\it singular vector} is a homogeneous nonzero vector $u$ in $M_\lmd$
such that
$$E_{i,i+1}(u)=0\qquad \for\;\;i=1,...,n-1.\eqno(2.18)$$
Here we have used the fact that all positive root vectors are generated by simple positive root vectors. The Verma module is irreducible if and only if any singular vecor is a scalar multiple of $v_\lmd$.

Consider the polynomial algebra
$${\cal A}=\mbb{C}[x_{i,j}\mid 1\leq j<i\leq n]\eqno(2.19)$$
in $n(n-1)/2$ variables. Set
$$x^\al=\prod_{1\leq j<i\leq n}x_{i,j}^{\al_{i,j}}\qquad\for\;\;\al\in\G.\eqno(2.20)$$
Then
$$\{x^\al\mid \al\in\G\}\;\;\mbox{forms a basis of}\;\;{\cal A}.\eqno(2.21)$$
Thus we have a linear isomorphism $\tau: M_{\lmd}\rta {\cal A}$ determined by
$$\tau(E^\al v_\lmd)=x^\al\qquad\for\;\;\al\in\G.\eqno(2.22)$$
The algebra ${\cal A}$ becomes $sl(n)$-module by the action
$$A(f)=\tau(A(\tau^{-1}(f)))\qquad\for\;\;A\in sl(n),\;f\in {\cal A}.\eqno(2.23)$$
For convenience, denote the partial derivatives
$$\ptl_{i,j}=\ptl_{x_{i,j}}\qquad\for\;\;1\leq j<i\leq n.\eqno(2.24)$$ 
In particular,
\begin{eqnarray*}& &d_i=E_{i,i+1}|_{\cal A}\\ &=&(\lmd_i-\sum_{j=i+1}^nx_{j,i}\ptl_{j,i}+\sum_{j=i+2}^nx_{j,i+1}\ptl_{j,i+1})\ptl_{i+1,i}+\sum_{j=1}^{i-1}x_{i,j}\ptl_{i+1,j}-\sum_{j=i+2}^nx_{j,i+1}\ptl_{j,i}\hspace{0.9cm}(2.25)\end{eqnarray*}
for $i=1,2,...,n-1$ by (2.15).
\psp

{\bf Proposition 2.1}. {\it A homogeneous vector} $u\in M_\lmd$ {\it is a singular vector if and only if} 
$$d_i(\tau(u))=0\qquad\for\;\;i=1,2,...,n-1.\eqno(2.26)$$
\pse

The system of partial differential equations
\begin{eqnarray*}\hspace{2.5cm}& &(\lmd_i-\sum_{j=i+1}^nx_{j,i}\ptl_{j,i}+\sum_{j=i+2}^nx_{j,i+1}\ptl_{j,i+1})\ptl_{i+1,i}(z)\\ & &+\sum_{j=1}^{i-1}x_{i,j}\ptl_{i+1,j}(z)-\sum_{j=i+2}^nx_{j,i+1}\ptl_{j,i}(z)=0\hspace{4.4cm}(2.27)\end{eqnarray*}
for $i=1,2,...,n-1$ and the unknown function $z$ in $\{x_{i,j}\mid 1\leq j<i\leq n\}$,  is called the {\it system of partial differential equations for singular vectors of type} $A_{n-1}$.

Next we want to find a family of exact power series solutions of the system
(2.27). First, we have
$$\eta_i=E_{i+1,i}|_{\cal A}=x_{i+1,i}+\sum_{j=1}^{i-1}x_{i+1,j}\ptl_{i,j}\eqno(2.28)$$
for $i=1,2,...,n-1$ by (2.16). Now we view $\{d_i,\eta_i\mid i=1,2,...,n-1\}$ purely as differential operators acting on functions of $\{x_{i,j}\mid 1\leq j<i\leq n\}$. In this way, we get a Lie algebra action on functions of $\{x_{i,j}\mid 1\leq j<i\leq n\}$ through $E_{i,i+1}=d_i$ and $E_{i+1,i}=\eta_i$ because
$sl(n)$ is generated by $\{E_{i,i+1},E_{i+1,i}\mid i=1,2,...,n-1\}$ as a Lie algebra.  Note that
$$h_i(E^\al v_\lmd)=(\lmd_i+\sum_{p=1}^{i-1}(\al_{i,p}-\al_{i+1,p})+\sum_{j=i+2}^n(\al_{j,i+1}-\al_{j,i})-2\al_{i+1,i})E^\al v_\lmd\eqno(2.29)$$
for $i=1,2,...,n-1$ and $\al\in\G$. Accordingly, we set
$$\zeta_i=h_i|_{\cal A}=\lmd_i+\sum_{p=1}^{i-1}(x_{i,p}\ptl_{i,p}-x_{i+1,p}\ptl_{i+1,p})+\sum_{j=i+2}^n(x_{j,i+1}\ptl_{j,i+1}-x_{j,i}\ptl_{j,i})-2x_{i+1,i}\ptl_{i+1,i}\eqno(2.30)$$
for $i=1,2,...,n-1$. The elements $h_i$ act on functions of $\{x_{i,j}\mid 1\leq j<i\leq n\}$ through $\zeta_i$.  A function $f$ of $\{x_{i,j}\mid 1\leq j<i\leq n\}$ is called {\it weighted} if there exist constants $\mu_1,\;\mu_2,\;...,\;\mu_{n-1}$ such that
$$\zeta_i(f)=\mu_if\qquad\for\;\;i=1,2,...,n-1.\eqno(2.31)$$
Since $d_i$ maps weigted functions to weigted functions, the system (2.27) is a weighted system. Any nonzero weighted solution of the system (2.27) is a singular vector of $sl(n)$. In particular, any nonzero  weighted  polynomial solution $f$ of the system (2.27) gives a singular vector $\tau^{-1}(f)$ in the Verma module $M_{\lmd}$.

Let
$${\cal A}_0=\mbb{C}[x_{i,j}\mid 1\leq j\leq i-2\leq n-2]\eqno(2.32)$$
be the polynomial algebra in $\{x_{i,j}\mid 1\leq j\leq i-2\leq n-2\}$. We denote
$$x^{\vec \mu}=\prod_{i=1}^{n-1}x_{i+1,i}^{\mu_i}\qquad\for\;\;\vec\mu=(\mu_1,\mu_2,...,\mu_{n-1})\in\mbb{C}^{\:n-1}.\eqno(2.33)$$
Let
$${\cal A}_1=\{\sum_{\vec j\in\mbb{N}^{\:n-1}}\sum_{i=1}^pf_{\vec\mu\:^i-\vec j}
x^{\vec\mu\:^i-\vec j}\mid 1\leq p\in\mbb{N},\;\vec \mu\:^i\in \mbb{C}^{\:n-1},\;f_{\vec\mu\:^i-\vec j}\in{\cal A}_0\}\eqno(2.34)$$
be the space of trucated-up formal power series in $\{x_{2,1},x_{3,2},...,x_{n,n-1}\}$ over ${\cal A}_0$. Then ${\cal A}$ is a subspace of ${\cal A}_1$. Since ${\cal A}_1$ is invariant under the action of $\{E_{i,i+1}=d_i, E_{i+1,i}=\eta_i\mid i=1,2,...,n-1\}$, ${\cal A}_1$ becomes an $sl(n)$-module.

For $\mu\in \mbb{C}$ and $p\in\mbb{N}$, we denote
$$\la \mu\ra_p=\mu(\mu-1)(\mu-2)\cdots (\mu-p+1).\eqno(2.35)$$
 Moreover,  by (2.28),  we define
$$\eta^{\mu}_i=(x_{i+1,i}+\sum_{j=1}^{i-1}x_{i+1,j}\ptl_{i,j})^{\mu}=\sum_{p=0}^{\infty}\frac{\la\mu\ra_p}{p!}x_{i+1,i}^{\mu-p}(\sum_{j=1}^{i-1}x_{i+1,j}\ptl_{i,j})^p\eqno(2.36)$$
as differential operators on ${\cal A}_1$, for $i=1,2,...,n-1$ and $\mu\in\mbb{C}$. If $\mu\not\in\mbb{N}$, then the above summation is infinite  and the positions of $x_{i+1,i}$ and $(\sum_{j=1}^{i-1}x_{i+1,j}\ptl_{i,j})$ are not symmetric. Since $x_{i+1,i}$ and $(\sum_{j=1}^{i-1}x_{i+1,j}\ptl_{i,j})$ commute, we have
$$\eta_i^{\mu_1}\eta_i^{\mu_2}=\eta^{\mu_1+\mu_2}_i\qquad\for\;\;\mu_1,\mu_2\in\mbb{C}.\eqno(2.37)$$

Given two differential operators $d$ and $\bar{d}$, we define the commutator
$$[d,\bar{d}]=d\bar{d}-\bar{d}d.\eqno(2.38)$$
 For any element $f\in {\cal A}_1$ and $r\in\mbb{C}$, we have
$$[\ptl_{i+1,i},x_{i+1,i}^r](f)=\ptl_{i+1,i}(x_{i+1,i}^rf)-x_{i+1,i}^r\ptl_{i+1,i}(f)=rx_{i+1,i}^{r-1}f,\eqno(2.39)$$
that is,
$$[\ptl_{i+1,i},x_{i+1,i}^r]=rx_{i+1,i}^{r-1}\qquad\mbox{as operators}.\eqno(2.40)$$
Note that if $(r,s)\not\in\{(i+1,j)\mid j=1,...,i\}$ and $(p,q)\not\in\{(i,j)\mid j=1,...,i-1\}$, then
$$[\ptl_{r,s},\eta^{\mu}_i]=[x_{p,q},\eta^{\mu}_i]=0\eqno(2.41)$$
directly by (2.36). Now
\begin{eqnarray*}\qquad[\ptl_{i+1,i},\eta^{\mu}_i]&=&\sum_{p=0}^{\infty}\frac{\la\mu\ra_p}{p!}[\ptl_{i+1,i},x_{i+1,i}^{\mu-p}(\sum_{j=1}^{i-1}x_{i+1,j}\ptl_{i,j})^p]\\ &=&\sum_{p=0}^{\infty}\frac{\la\mu\ra_p}{p!}(\mu-p)x_{i+1,i}^{\mu-p-1}(\sum_{j=1}^{i-1}x_{i+1,j}\ptl_{i,j})^p\\ &=&\sum_{p=0}^{\infty}\frac{\mu \la\mu-1\ra_p}{p!}x_{i+1,i}^{\mu-p-1}(\sum_{j=1}^{i-1}x_{i+1,j}\ptl_{i,j})^p
=\mu \eta^{\mu-1}_i\hspace{3.3cm}(2.42)\end{eqnarray*}
by (2.40). Moreover, for $j=1,2,...,i-1$,
\begin{eqnarray*}\qquad [\ptl_{i+1,j},\eta^{\mu}_i] &=&\sum_{p=0}^{\infty}\frac{\la\mu\ra_p}{p!}[\ptl_{i+1,j},x_{i+1,i}^{\mu-p}(\sum_{j=1}^{i-1}x_{i+1,j}\ptl_{i,j})^p]\\ &=& \sum_{p=0}^{\infty}\frac{\la\mu\ra_p}{p!}px_{i+1,i}^{\mu-p}(\sum_{j=1}^{i-1}x_{i+1,j}\ptl_{i,j})^{p-1}\ptl_{i,j}\\ &=&\sum_{p=1}^{\infty}\frac{\mu\la\mu-1\ra_{p-1}}{(p-1)!}x_{i+1,i}^{\mu-p}(\sum_{j=1}^{i-1}x_{i+1,j}\ptl_{i,j})^{p-1}\ptl_{i,j} =\mu \eta^{\mu-1}_i\ptl_{i,j}\hspace{1.8cm}(2.43)\end{eqnarray*}
and similarly,
$$[x_{i,j},\eta^{\mu}_i]=-\mu \eta^{\mu-1}_ix_{i+1,j}\qquad\for\;\;j=1,2,...,i-1.\eqno(2.44)$$

{\bf Lemma 2.2}. {\it For} $i,l\in\{1,2,...,n-1\}$ {\it and} $\mu\in\mbb{C}$, {\it we have}:
$$[d_l,\eta_i^{\mu}]=\mu\dlt_{i,l} \eta_i^{\mu-1}(1-\mu+\zeta_i).\eqno(2.45)$$

{\it Proof}. Note that
$$[E_{l,l+1},E_{i+1,i}^m]=m\dlt_{i,l}E_{i+1,i}^{m-1}(1-m+h_i)\qquad\for\;\;m\in\mbb{N}\eqno(2.46)$$
(cf. (2.3)). So (2.45) holds for any $\mu\in\mbb{N}$ by (2.25), (2.28) and (2.30). Since (2.45) is completely determined by  (2.41)-(2.44), which are independent of whether $\mu$ is a nonnegative integer, it must hold for any $\mu\in\mbb{C}.\qquad\Box$
\psp

Denote the  Cartan matrix of $sl(n)$ by
$$\left[\begin{array}{cccc}a_{1,1}&a_{1,2}&\cdots& a_{1,{n-1}}\\ a_{2,1}&a_{2,2}&\cdots &a_{2,{n-1}}\\ \vdots&\vdots& &\vdots\\ a_{n-1,1}&a_{n-1,2}&\cdots& a_{n-1,{n-1}}\end{array}\right]=\left[\begin{array}{rrrr}2&-1&&\\ -1&2&\ddots &\\ &\ddots& \ddots&-1\\ &&-1&2\end{array}\right].\eqno(2.47)$$

{\bf Lemma 2.3}. {\it For} $i,l\in\{1,2,...,n-1\}$ {\it and} $\mu\in\mbb{C}$, 
$$[\zeta_l,\eta_i^{\mu}]=-\mu a_{l,i}\eta_i^{\mu}.\eqno(2.48)$$

{\it Proof}. Observe that
$$[h_l,E_{i+1,i}^m]=-ma_{l,i}E_{i+1,i}^m\qquad\for\;\;m\in\mbb{N}\eqno(2.49)$$
(cf. (2.3)). Hence (2.48) holds for any $\mu\in\mbb{N}$ by  (2.28) and (2.30). Again  (2.48) is completely determined by  (2.41)-(2.44), which are independent of whether $\mu$ is a nonnegative integer. Thus (2.48) must hold for any $\mu\in\mbb{C}.\qquad\Box$
\psp

Let $m$ be a positive integer and let 
$$\Im: \{1,2,...,m\}\rta\{1,2,...,n-1\}\eqno(2.50)$$
be a map such that
$$\Im(i)\neq \Im(i+1).\eqno(2.51)$$
We define
$$\iota_1=\lmd_{\Im(1)}+1\eqno(2.52)$$
and
$$\iota_i=\lmd_{\Im(i)}+1-\sum_{p=1}^{i-1}a_{\Im(i),\Im(p)}\iota_p\eqno(2.53)$$
for $i>1$ by induction on $i$. Set
$$\eta[\Im]=\eta_{\Im(m)}^{\iota_m}\eta_{\Im(m-1)}^{\iota_{m-1}}\cdots \eta_{\Im(2)}^{\iota_2}\eta_{\Im(1)}^{\iota_1}(1),\eqno(2.54)$$
a product of differential operators acting on $1$. Then
$$\eta[\Im]\in {\cal A}_1\eqno(2.55)$$
(cf. (2.34)). Fot instance,
$$\eta_{\Im(1)}^{\iota_1}(1)=x_{\Im(1)+1,\Im(1)}^{\lmd_{\Im(1)}+1}.\eqno(2.56)$$
By Lemma 2.2 and Lemma 2.3, we have:
\psp

{\bf Theorem 2.4}. {\it The power series} $\eta[\Im]$ {\it is a weighted exact solution of the system (2.27).}

\section{Complete Solutions}

In this section, we want to solve (2.27) in ${\cal A}_1$ completely. 

For any
$$\vec\mu=(\mu_1,\mu_2,...,\mu_{n-1})\in \mbb{C}^{\:n-1},\eqno(3.1)$$
\newpage

we define
\begin{eqnarray*}\hspace{1.2cm}\phi_{\vec \mu}&=&\eta_2^{\mu_2}\eta_1^{\mu_2-\lmd_{n-1}-1}\cdots \eta_i^{\mu_i}\eta_{i-1}^{\mu_i-\lmd_{n-1}-1}\cdots\eta_1^{\mu_i-\lmd_{n-1}-\cdots-\lmd_{n-(i-1)}-(i-1)}\\ & &\cdots \eta_{n-1}^{\mu_{n-1}}\eta_{n-2}^{\mu_{n-1}-\lmd_{n-1}-1}\cdots\eta_1^{\mu_{n-1}+2-\lmd_2-\cdots-\lmd_{n-1}-n}(1).\hspace{3.7cm}(3.2)\end{eqnarray*}
Then
$$\phi_{\vec \nu}\in {\cal A}_1\eqno(3.3)$$
(cf. (2.32)-(2.34)) and
is a solution of the system
$$d_2(z)=d_3(z)=\cdots=d_{n-1}(z)=0\eqno(3.4)$$
by Lemma 2.2 and 2.3. In fact, we have:
\psp

{\bf Lemma 3.1}. {\it An element} $z$ {\it  in} ${\cal A}_1$ {\it is a solution of the system (3.4) if and only if it can be written as}
$$z=\sum_{\vec j\in\mbb{N}^{\:n-1}}\sum_{i=1}^pc_{\vec\mu\:^i-\vec j}x_{2,1}^{\mu^i_1-j_1}\phi_{\vec \mu\:^i-\vec j}\qquad\mbox{\it  with}\;\;c_{\vec\mu-\:^i\vec j}\in \mbb{C}\eqno(3.5)$$
{\it for some} $\vec\mu\:^1,...,\vec\mu\:^p\in \mbb{C}^{\:n-1}$.
\pse

{\it Proof}. By Lemma 2.2 and 2.3, we only need to prove the necessity. Recall that 
$$E_{i,i+1}=d_i,\;\;E_{i+1,i}=\eta_i\qquad\mbox{as operators on}\;\;{\cal A}_1\eqno(3.6)$$
(cf. (2.25) and (2.28)) for $i=1,2,...,n-1$. Note
$$d_{n-1}=(\lmd_{n-1}-x_{n,n-1}\ptl_{n,n-1})\ptl_{n,n-1}+\sum_{i=1}^{n-2}x_{n-1,i}\ptl_{n,i}.\eqno(3.7)$$
Moreover, (2.11) tells us that
\begin{eqnarray*}\hspace{1.5cm}& &d_{n-2,n}=E_{n-2,n}|_{{\cal A}_1}\\ &=&
(\lmd_{n-1}+\lmd_{n-2}-x_{n,n-2}\ptl_{n,n-2}-x_{n,n-1}\ptl_{n,n-1})\ptl_{n,n-2}\\ & &-d_{n-1}\ptl_{n-1,n-2}+\sum_{i=1}^{n-3}x_{n-2,i}\ptl_{n,i}\\ &=&(\lmd_{n-1}+\lmd_{n-2}+1-x_{n,n-2}\ptl_{n,n-2}-x_{n,n-1}\ptl_{n,n-1})\ptl_{n,n-2}\\ & &-\ptl_{n-1,n-2}d_{n-1}+\sum_{i=1}^{n-3}x_{n-2,i}\ptl_{n,i}.\hspace{7.2cm}(3.8)\end{eqnarray*}
Set
$$\bar{\lmd}_i=n-i-1+\sum_{p=i}^{n-1}\lmd_p\qquad\for\;\;i=2,3,...,n-2.\eqno(3.9)$$
Furthermore,
\begin{eqnarray*}\qquad& &d_{n-3,n}=E_{n-3,n}|_{{\cal A}_1}\\ &=& (\bar{\lmd}_{n-3}-2-x_{n,n-3}\ptl_{n,n-3}-x_{n,n-2}\ptl_{n,n-2}-x_{n,n-1}\ptl_{n,n-1})\ptl_{n,n-3}\\ & &-d_{n-2,n}\ptl_{n-2,n-3}-d_{n-1}\ptl_{n-1,n-3}+\sum_{i=1}^{n-4}x_{n-3,i}\ptl_{n,i}\\ &=&(\bar{\lmd}_{n-3}-2-\sum_{p=1}^3x_{n,n-p}\ptl_{n,n-p}))\ptl_{n,n-3} -\ptl_{n-2,n-3}d_{n-2,n}-\ptl_{n-1,n-3}d_{n-1,n}\\& &-[d_{n-2,n},\ptl_{n-2,n-3}]-[d_{n-1},\ptl_{n-1,n-3}]+\sum_{i=1}^{n-4}x_{n-3,i}\ptl_{n,i}\\ &=&
(\bar{\lmd}_{n-3}-\sum_{p=1}^3x_{n,n-p}\ptl_{n,n-p})\ptl_{n,n-3}+\sum_{i=1}^{n-4}x_{n-3,i}\ptl_{n,i}\\ & &-\ptl_{n-2,n-3}d_{n-2,n}-\ptl_{n-1,n-3}d_{n-1,n}.\hspace{7.1cm}(3.10)\end{eqnarray*}
By induction, we can prove that
$$d_{i,n}=E_{i,n}|_{{\cal A}_1}=(\bar{\lmd}_i-\sum_{p=i}^{n-1}x_{n,p}\ptl_{n,p})\ptl_{n,i}+\sum_{q=1}^{i-1}x_{i,q}\ptl_{n,q}-\sum_{j=i+1}^{n-1}\ptl_{j,i}d_{j,n}\eqno(3.11)$$
for $i=2,3,...,n-2$, where we take
$$d_{n-1,n}=d_{n-1}.\eqno(3.12)$$

Suppose that  
$$z=fx_{2,1}^{\mu_1}\phi_{\vec \mu}\eqno(3.13)$$
is a solution of the system (3.4) for some $\vec\mu\in \mbb{C}^{\:n-1}$ and $f\in{\cal A}_0$ (cf. (2.32) and (3.2)). We want to prove that  $f$ is a constant.  Denote
$$\es_i=(0,...,0,\stl{i}{1},0,...,0)\in\mbb{C}^{\:n-1}\eqno(3.14)$$
and 
$$\iota_{i,j}=\mu_{i+j}-\sum_{p=1}^j(\lmd_{n-p}+1)\qquad\for\;\;2\leq i\leq n-1,\;0\leq j\leq n-i-1.\eqno(3.15)$$
We define
$$U_i=\{\sum_{\vec j\in\mbb{N}^{\:n-1}}g_{\vec j}x_{2,1}^{\mu_1-j_1}\phi_{\vec \mu-\es_i-\vec j}\mid g_{\vec j}\in {\cal A}_0\}\eqno(3.16)$$
for $i=2,...,n-1$, and
$$U=\sum_{i=2}^{n-1}U_i.\eqno(3.17)$$
Moreover, for fixed $i\geq 2$,  
$$\{\eta_i^{\iota_{i,j}}\mid 0\leq j\leq n-i-1\}\eqno(3.18)$$
are all the factors in the righthand side of (3.2) that contain $x_{i+1,p}$ or
$\ptl_{i,q}$ with $p=1,...,i$ and $q=1,2,....,i-1$ by (2.36). Besides,
$$[\ptl_{i+1,i},\eta_i^{\iota_{i,j}}]=\iota_{i,j}\eta_i^{\iota_{i,j}-1},\eqno(3.19)$$
$$[\ptl_{i+1,p},\eta_i^{\iota_{i,j}}]=\iota_{i,j}\eta_i^{\iota_{i,j}-1}\ptl_{i,p}\qquad\for\;\;p=1,...,i-1,\eqno(3.20)$$
$$[\sum_{p=r}^ix_{i+1,p}\ptl_{i+1,p},\eta_i^{\iota_{i,j}}]=\iota_{i,j}(\eta_i^{\iota_{i,j}}-\sum_{p=1}^{r-1}x_{i+1,p}\eta_i^{\iota_{i,j}-1}\ptl_{i,p}),\eqno(3.21)$$
$$[\sum_{p=q}^{i-1}x_{i,p}\ptl_{i,p},\eta_i^{\iota_{i,j}}]=-\iota_{i,j}\sum_{p=q}^{i-1}x_{i+1,p}\eta_i^{\iota_{i,j}-1}\ptl_{i,p},\eqno(3.22)$$
where $1\leq r\leq i$ and $1\leq q\leq i-1$.

By (3.18)-(3.22), we have
$$\ptl_{i+1,r}(\phi_{\vec \mu})\in  U_i\eqno(3.23)$$
and
$$(\sum_{p=r}^ix_{i+1,p}\ptl_{i+1,p})(\phi_{\vec \mu})\equiv c_{i,r}\phi_{\vec \mu}
\;\;({\rm mod}\;\sum_{s=1}^iU_s)\eqno(3.24)$$
for $2\leq i\leq n-1$ and $1\leq r\leq i$. Since for $2\leq i\leq n-2$,
$$E_{i,n}=[E_{i,i+1},[E_{i+1,i+2},\cdots,[E_{n-2,n-1},E_{n-1,n}]\cdots]],\eqno(3.25)$$
we have
$$d_{i,n}=[d_i,[d_{i+1},\cdots [d_{n-2},d_{n-1}]\cdots]].\eqno(3.26)$$
Thus
$$d_{i,n}(z)=0\qquad\for\;\;2\leq i\leq n-1.\eqno(3.27)$$
By (3.11), (3.16) and (3.22)-(3.27), we obtain
$$d_{i,n}(z)=[(\bar{\lmd}_i-c_{n-1,i})\ptl_{n,i}(f)+\sum_{q=1}^{i-1}x_{i,q}\ptl_{n,q}(f)]x_{2,1}^{\mu_1}\phi_{\vec \mu}-\sum_{j=i+1}^{n-1}\ptl_{j,i}d_{j,n}(z)\equiv 0\;\;({\rm mod}\;U)\eqno(3.28)$$
for  $i=2,3,...n-2$ and
$$d_{n-1}(z)=(\sum_{q=1}^{n-2}x_{n-1,q}\ptl_{n,q}(f))x_{2,1}^{\mu_1}\phi_{\vec \mu}\equiv 0\;\;({\rm mod}\;U).\eqno(3.29)$$
 Since the constraint on $d_{i,n}(z)\equiv 0\;\;({\rm mod}\;U)$ implies $\ptl_{i,n}d_{i,n}(z)\equiv 0\;\;({\rm mod}\;U)$, (3.28) is equivalent to
$$d_{i,n}(z)=[(\bar{\lmd}_i-c_{n-1,i})\ptl_{n,i}(f)+\sum_{q=1}^{i-1}x_{i,q}\ptl_{n,q}(f)]x_{2,1}^{\mu_1}\phi_{\vec \mu}\equiv 0\;\;({\rm mod}\;U)\eqno(3.30)$$
for  $i=2,3,...n-2$.

Expressions (3.29) and (3.30) give
$$\sum_{q=1}^{i-1}x_{i,q}\ptl_{n,q}(f)+(\bar{\lmd}_i-c_{n-1,i})\ptl_{n,i}(f)=0\qquad\for\;\;i=2,...,n-2,\eqno(3.31)$$
$$\sum_{q=1}^{n-2}x_{n-1,q}\ptl_{n,q}(f)=0.\eqno(3.32)$$
We view 
$$\ptl_{n,1}(f),\;\ptl_{n,2}(f),\;\cdots,\;\ptl_{n,n-2}(f)\;\;\mbox{ as unknowns}.\eqno(3.33)$$
Then the coefficient determinant of the system (3.31) and (3.32) is
\begin{eqnarray*}\hspace{2cm}& &\left|\begin{array}{cccc}x_{2,1}&\bar{\lmd}_2-c_{n-1,2}& &\\ x_{3,1}& x_{3,2}&\ddots&\\ \vdots&\ddots&\ddots &\bar{\lmd}_{n-2}-c_{n-1,n-2}\\ x_{n-1,1}&\cdots &x_{n-1,n-3}& x_{n-1,n-2}\end{array}\right|\\ &=&\prod_{p=2}^{n-1}x_{p,p-1}-g(x_{2,1},x_{3,2},...,x_{n-1,n-2})\not\equiv 0,\hspace{4.7cm}(3.34)\end{eqnarray*}
where $g(x_{2,1},x_{3,2},...,x_{n-1,n-2})$ is a polynomial of degree $n-3$ in  $\{x_{2,1},x_{3,2},...,x_{n-1,n-2}\}$ over ${\cal A}_0$ (cf. (2.32)). 
Therefore,
$$\ptl_{n,q}(f)=0\qquad\for\;\;q=1,2,....,n-2.\eqno(3.35)$$

Based on our calculations in (3.23)-(3.25), we can prove by induction that
$$\ptl_{q+r,q}(f)=0\qquad\for\;\;1\leq q\leq n-2,\;2\leq r\leq n-q.\eqno(3.36)$$
So  $f$ is  a constant.

Suppose that $z$ is any solution of the system (3.4) in ${\cal A}_1$. By (2.34) and (3.2), it can be written as 
$$z=\sum_{\vec j\in\mbb{N}^{\:n-1}}\sum_{i=1}^pf_{\vec\mu^i-\vec j}
x_{2,1}^{\mu^i_1-j_1}\phi_{\vec \mu^i-\vec j}\qquad\mbox{ with}\;\;f_{\vec j}\in {\cal A}_0\eqno(3.37)$$
Let 
$$S=\{\vec\sgm\in\mbb{C}^{\:n-1}\mid f_{\vec\sgm}\neq 0;\;f_{\vec\sgm+\vec j}=0\;\mbox{for all}\;\vec 0\neq\vec j\in\mbb{N}^{\:n-1}\}.\eqno(3.38)$$
The above arguments show that
$$\{f_{\vec \sgm}\mid \vec\sgm\in S\}\qquad \mbox{are constants}\eqno(3.39)$$
(cf. the key equations (3.29) and (3.30)). Since
$\sum_{\sgm\in S}f_{\sgm}x_{2,1}^{\sgm_1}\phi_{\sgm}$
is a solution of the system (3.4), so is $z-\sum_{\sgm\in S}f_{\sgm}x_{2,1}^{\sgm_1}\phi_{\sgm}$. By induction, we prove the lemma.$\qquad\Box$
\psp

To solve the system (2.27) in ${\cal A}_1$, we only need to consider the solutions of the form $z=x_{2,1}^{\mu_1}\phi_{\vec\mu}$ with $\vec\mu\in \mbb{C}^{\:n-1}$ by the above lemma because (2.27) is a weighted system. Note
$$d_1=(\lmd_1-\sum_{j=2}^nx_{j,1}\ptl_{j,1}+\sum_{j=3}^nx_{j,2}\ptl_{j,2})\ptl_{2,1}-\sum_{j=3}^nx_{j,2}\ptl_{j,1}.\eqno(3.40)$$
Denote
$$\mu_{1,r}=\sum_{p=r}^{n-1}\mu_i-\sum_{p=2}^{n-r}(p-1)(\lmd_p+1)-(n-r)\sum_{q=n-r+1}^{n-1}(\lmd_q+1),\;\;\mu_{1,1}=\mu_1+\mu_{1,2}\eqno(3.41)$$
for $r=2,3,...,n-1$,
$$\mu_{2,r}=\mu_r-\sum_{p=1}^{r-2}(\lmd_{n-p}+1)\qquad \for\;\;r=2,3,...,n-1\eqno(3.42)$$
and
$$\td{\mu}=\sum_{i=2}^{n-1}\mu_i-\sum_{p=3}^{n-1}(p-2)(\lmd_p+1).\eqno(3.43)$$
Letting $x_{p,q}=0$ for $1\leq q\leq p-2\leq n-2$ in
$$d_1(z)=[(\lmd_1-\sum_{j=2}^nx_{j,1}\ptl_{j,1}+\sum_{j=3}^nx_{j,2}\ptl_{j,2})\ptl_{2,1}-\sum_{j=3}^nx_{j,2}\ptl_{j,1}](x_{2,1}^{\mu_1}\phi_{\vec\mu})=0,\eqno(3.44)$$
we get
$$\mu_{1,1}(\lmd_1+1-\mu_{1,1}+\td{\mu})-\sum_{r=2}^{n-1}\mu_{2,r}\mu_{1,r}=0\eqno(3.45)$$
by (2.36) and (3.2). 

Suppose  $n>3$. We take $x_{p,q}=0$ for $1\leq q\leq p-2\leq n-2$ in
$$\ptl_{n,1}d_1(z)=\ptl_{n,1}[(\lmd_1-\sum_{j=2}^nx_{j,1}\ptl_{j,1}+\sum_{j=3}^nx_{j,2}\ptl_{j,2})\ptl_{2,1}-\sum_{j=3}^nx_{j,2}\ptl_{j,1}](x_{2,1}^{\mu_1}\phi_{\vec\mu})=0,\eqno(3.46)$$ 
and obtain
\begin{eqnarray*} \hspace{2cm}& &\mu_{n-1}\left[\prod_{i=1}^{n-2}(\mu_{n-1}-i-\sum_{p=1}^i\lmd_{n-p})\right][(\mu_{1,1}-1)(\lmd_1-\mu_{1,1}+\td{\mu})\\ & &-\sum_{r=2}^{n-2}\mu_{2,r}(\mu_{1,r}-1)-(\mu_{2,n-1}-1)(\mu_{1,n-1}-1)]=0.\hspace{3.2cm}(3.47)\end{eqnarray*}
Note that
\begin{eqnarray*}& &[\mu_{1,1}(\lmd_1+1-\mu_{1,1}+\td{\mu})-\sum_{r=2}^{n-1}\mu_{2,r}\mu_{1,r}]-[(\mu_{1,1}-1)(\lmd_1-\mu_{1,1}+\td{\mu})\\ &&-\sum_{r=2}^{n-2}\mu_{2,r}(\mu_{1,r}-1)-(\mu_{2,n-1}-1)(\mu_{1,n-1}-1)]\hspace{8cm}\end{eqnarray*}
\begin{eqnarray*} &=&\lmd_1+\td{\mu}-\sum_{r=2}^{n-1}\mu_{2,r}+1-\mu_{1,n-1}=\lmd_1+1-\mu_{1,n-1}\\ &=& (n-1)+\sum_{i=1}^{n-1}\lmd_i-\mu_{n-1}.\hspace{9.6cm}(3.48)\end{eqnarray*}
By (3.45), (3.47) and (3.48), we have
$$\mu_{n-1}\prod_{i=1}^{n-1}(\mu_{n-1}-i-\sum_{p=1}^i\lmd_{n-p})=0.\eqno(3.49)$$
Therefore,
$$\mu_{n-1}\in\{0,i+\sum_{p=1}^i\lmd_{n-p}\mid i=1,2,...,n-1\}.\eqno(3.50)$$

Assume $n=3$. Then 
$$d_1=(\lmd_1-x_{2,1}\ptl_{2,1}-x_{3,1}\ptl_{3,1}+x_{3,2}\ptl_{3,2})\ptl_{2,1}-x_{3,2}\ptl_{3,1},\eqno(3.51)$$
$$z=x_{2,1}^{\mu_1}\phi_{\vec\mu}=x_{2,1}^{\mu_1}(x_{3,2}+x_{3,1}\ptl_{2,1})^{\mu_2}(x_{2,1}^{\mu_2-\lmd_2-1}),\eqno(3.52)$$
and (3.45) becomes
$$(\mu_1+\mu_2-\lmd_2-1)(\lmd_1+\lmd_2+2-\mu_1)-\mu_2(\mu_2-\lmd_2-1)=0\eqno(3.53)$$
Letting $x_{3,1}=0$ in
\begin{eqnarray*}&&\ptl_{3,1}d_1(z)=\ptl_{3,1}[(\lmd_1-x_{2,1}\ptl_{2,1}-x_{3,1}\ptl_{3,1}+x_{3,2}\ptl_{3,2})\ptl_{2,1}\\& &-x_{3,2}\ptl_{3,1}][x_{2,1}^{\mu_1}(x_{3,2}+x_{3,1}\ptl_{2,1})^{\mu_2}(x_{2,1}^{\mu_2-\lmd_2-1})]=0,\hspace{5.7cm}(3.54)\end{eqnarray*}
we get
\begin{eqnarray*}\hspace{1cm}& &\mu_2(\mu_2-\lmd_2-1)(\mu_1+\mu_2-\lmd_2-2)(\lmd_1+\lmd_2+1-\mu_1)\\ & &-\mu_2(\mu_2-1)(\mu_2-\lmd_2-1)(\mu_2-\lmd_2-2)=0,\hspace{5.2cm}(3.55)\end{eqnarray*}
equivalently
$$\mu_2(\mu_2-\lmd_2-1)[(\mu_1+\mu_2-\lmd_2-2)(\lmd_1+\lmd_2+1-\mu_1)-(\mu_2-1)(\mu_2-\lmd_2-2)]=0.\eqno(3.56)$$
By (3.53), we have
\begin{eqnarray*}\hspace{1cm}& &(\mu_1+\mu_2-\lmd_2-2)(\lmd_1+\lmd_2+1-\mu_1)-(\mu_2-1)(\mu_2-\lmd_2-2)\\ &=& -(\lmd_1+\mu_2)+\mu_2(\mu_2-\lmd_2-1)
-(\mu_2-1)(\mu_2-\lmd_2-2)\\ &=&-(\lmd_1+\mu_2)+2\mu_2-\lmd_2-2\\ &=&\mu_2-\lmd_1-\lmd_2-2.\hspace{9.7cm}(3.57)\end{eqnarray*}
Thus (3.56) and (3.57) give
$$\mu_2(\mu_2-\lmd_2-1)(\mu_2-\lmd_1-\lmd_2-2)=0,\eqno(3.58)$$
which implies that (3.50) holds for any $n\geq 2$.

When $n=2$, the solution space of (2.27) is $\mbb{C}+\mbb{C}x_{2,1}^{\lmd_1+1}$. In general, we can use (3.50) to reduce the problem of solving (2.27) to $sl(n-1)$ as follows.  Denote
\begin{eqnarray*}\hspace{1.5cm}& &\Psi_i=x_{2,1}^{\mu_1}\eta_2^{\mu_2}\eta_1^{\mu_2-\lmd_{n-1}-1}\cdots \eta_i^{\mu_i}\eta_{i-1}^{\mu_i-\lmd_{n-1}-1}\cdots\eta_1^{\mu_i-\lmd_{n-1}-\cdots-\lmd_{n-(i-1)}-(i-1)}\\ & &\cdots \eta_{n-2}^{\mu_{n-2}}\eta_{n-3}^{\mu_{n-2}-\lmd_{n-1}-1}\cdots\eta_1^{\mu_{n-2}+3-\lmd_3-\cdots-\lmd_{n-1}-n}\\ & &\eta_{i-2}^{-\lmd_{i-1}-1}\eta_{i-3}^{-\lmd_{i-1}-\lmd_{i-2}-2}\cdots \eta_1^{-\lmd_2-\cdots-\lmd_{i-1}-(i-2)}\hspace{5.2cm}(3.59)\end{eqnarray*}
for $i=0,1,...,n-1$, where we treat
$$\eta_{i-2}^{-\lmd_{i-1}-1}\eta_{i-3}^{-\lmd_{i-1}-\lmd_{i-2}-2}\cdots \eta_1^{-\lmd_2-\cdots-\lmd_{i-1}-(i-2)}=1\qquad\mbox{if}\;\;i=0,1,2.\eqno(3.60)$$
Moreover, we set
$$\psi_0=1,\;\;\psi_i=\eta_{n-1}^{n-i+\sum_{p=1}^{n-i}\lmd_{n-p}}\eta_{n-2}^{n-i-1+\sum_{p=2}^{n-i}\lmd_{n-p}}\cdots\eta_i^{\lmd_i+1}(1),\eqno(3.61)$$
for $i=1,2,...,n-1$. Then $\{\psi_i\mid i=0,1,...,n-1\}$ are solutions of (2.27) by Theorem 2.4.  Denote
$$\lmd_{n-1,0}=0,\;\;\lmd_{n-1,i}=n-i+\sum_{p=1}^{n-i}\lmd_{n-p}\qquad\for\;\;i=1,2,...,n-1.\eqno(3.62)$$
According to (3.50),
$$\mu_{n-1}=\lmd_{n_1,i_{n-1}}\qquad\mbox{for some}\;\;i_{n-1}\in\{0,1,...,n-1\}.\eqno(3.63)$$
Thus,
$$z=x_{2,1}^{\mu_1}\phi_{\vec\mu}=\Psi_{i_{n-1}}\psi_{i_{n-1}}.\eqno(3.64)$$

Set
$$\lmd_j^{(n-2)}=\left\{\begin{array}{ll}\lmd_j&\mbox{if}\;\;j<i_{n-1}-1,\\ \lmd_{i_{n-1}}+\lmd_{i_{n-1}-1}+1&\mbox{if}\;\;j=i_{n-1}-1,\\ \lmd_{j+1}&\mbox{if}\;\;j\geq i_{n-1}\end{array}\right.\eqno(3.65)$$
for $j=1,2,...,n-2$. By (2.29) and (2.30),
$$h_j(\psi_{i_{n-1}})=\zeta_j (\psi_{i_{n-1}})=\lmd^{(n-2)}_j\psi_{i_{n-1}}\qquad\for\;\; j=1,2,...,n-2.\eqno(3.66)$$
Define
\begin{eqnarray*}\hspace{1cm}& &\bar{d_i}=(\lmd^{(n-1)}_i-\sum_{j=i+1}^{n-1}x_{j,i}\ptl_{j,i}+\sum_{j=i+2}^{n-1}x_{j,i+1}\ptl_{j,i+1})\ptl_{i+1,i}\\ &&+\sum_{j=1}^{i-1}x_{i,j}\ptl_{i+1,j}-\sum_{j=i+2}^{n-1}x_{j,i+1}\ptl_{j,i}\hspace{7.7cm}(3.67)\end{eqnarray*}
for $i=1,2,...,n-2$ by (2.25). Then
$$\bar{d}_i(\Psi_{i_{n-1}}(1))=0\qquad \for\;\;i=1,2,...,n-2\eqno(3.68)$$
by (3.64), (3.66) and the fact that $\psi_{i_{n-1}}$ is a  solution of (2.27). The system (3.68) is a version (2.27) for $sl(n-1)$. So we have reduced the problem of solving (2.27) for $sl(n)$ to that for $sl(n-1)$. This gives us an inductive process of solving (2.27) completely.

Suppose that $\{\lmd_j^{(p)}\mid j=1,2,,,,,p\}$ are defined. We define
$$\lmd_{p,0}=0,\;\;\lmd_{p,j}=p-j+1+\sum_{r=1}^{p-j+1}\lmd_{p-r+1}^{(p)}\qquad\for\;\;j=1,2,...,p.\eqno(3.69)$$
We choose 
$$i_p\in\{0,1,...,p\}\eqno(3.70)$$
and define
$$\lmd_j^{(p-1)}=\left\{\begin{array}{ll}\lmd_j^{(p)}&\mbox{if}\;\;j<i_p-1,\\ \lmd^{(p)}_{i_p}+\lmd^{(p)}_{i_p-1}+1&\mbox{if}\;\;j=i_p-1,\\ \lmd^{(p)}_{j+1}&\mbox{if}\;\;j\geq i_p\end{array}\right.\eqno(3.71)$$
By induction, we have defined all
$$\{\lmd^{(p)}_j,\lmd_{p,r}\mid p=1,2,...,n-1,\;j=1,2,...,p,\;r=0,1,...,p\}\eqno(3.72)$$
for a given vector
$$\vec i=(i_1,i_2,...,i_{n-1})\in\mbb{N}^{\:n-1}\;\;\mbox{with}\;\; i_p\leq p,\eqno(3.73)$$
where
$$\lmd^{(n-1)}_j=\lmd_j\qquad \for\;\;j=1,2,...,n-1.\eqno(3.74)$$
For $\vec i$ in (3.73), we set
\begin{eqnarray*}\hspace{2cm}& &\theta_{\vec i}=\eta_1^{\lmd_{1,i_1}}\cdots\eta_p^{\lmd_{p,i_p}}\eta_{p-1}^{\lmd_{p,i_p}-\lmd^{(p)}-1}\cdots \eta_{i_p}^{\lmd^{(p)}_{i_p}+1}\\ & &\cdots \eta_{n-1}^{\lmd_{n-1,i_{n-1}}}\eta_{n-2}^{\lmd_{n-1,i_{n-1}}-\lmd_{n-1}-1}\cdots \eta_{i_{n-1}}^{\lmd_{i_{n-1}}+1}(1),\hspace{4.4cm}(3.75)\end{eqnarray*}
where we treat
$$\eta_p^{\lmd_{p,i_p}}\eta_{p-1}^{\lmd_{p,i_p}-\lmd^{(p)}-1}\cdots \eta_{i_p}^{\lmd^{(p)}_{i_p}+1}=1\qquad\mbox{if}\;\;i_p=0.\eqno(3.76)$$
Note that $\theta_{\vec i}$ is of the form $\eta[\Im]$ defined in (2.50)-(2.54). In fact, all exponent of any $\eta_j$ in (2.54) is either 0 or of the form
$\sum_{r=1}^s(\lmd_{q+r}+1)$ by (3.62), (3.65), (3.69), (3.71) and (3.75). The following is our main theorem in this section.
\psp

{\bf Theorem 3.2}. {\it The solution space the system (2.27) is}:
$${\cal S}=\sum_{\vec i=(i_1,...,i_{n-1})\in\mbb{N}^{\:n-1},\;i_p\leq p}\mbb{C}\theta_{\vec i},\eqno(3.77)$$
{\it which has dimension} $n!$.
\psp

Assume that
$$\ves=n-2+\sum_{p=1}^{n-1}\lmd_p\in\mbb{N}.\eqno(3.78)$$
Set
$$\ves_q=q+\sum_{p=1}^q\lmd_p\qquad\for\;\;q=1,2,...,n-2.\eqno(3.79)$$
Then
$$\lmd_{n-1,q+1}+\ves_q=\ves\qquad\for\;\;q=1,2,...,n-2\eqno(3.80)$$
by (3.62). Define
$$\phi=\eta_1^{\lmd_{n-1,2}}\eta_2^{\lmd_{n-1,3}}\cdots\eta_{n-2}^{\lmd_{n-1,n-1}}\eta_{n-1}^{\ves}\eta_{n-2}^{\ves_{n-2}}\cdots \eta_2^{\ves_2}\eta_1^{\ves_1}(1),\eqno(3.81)$$
which is a special case of $\theta_{\vec i}$. It can be verified that
$\phi$ is a polynomial. Thus $\tau^{-1}(\phi)$ (cf. (2.22)) is a nontrivial singular vector in the Verma module $M_\lmd$ (cf. (2.14)), which was obtained by
Malikov, Feigin and Fuchs [MFF].

Suppose that $\theta_{\vec i}$ is a polynomial. If 
$$\lmd_{n-1,j}=n-j+\sum_{p=1}^{n-j}\lmd_{n-p}\not\in \mbb{N}+1\qquad\for\;\;j
=1,2,...,n-1,\eqno(3.82)$$
then $i_{n-1}=0$. By induction, we obtain:
\psp

{\bf Corollary 3.3}. {\it The Verma module} $M_\lmd$ {\it is irreducible if and only if}
$$j+\sum_{p=1}^{j}\lmd_{i+p}\not\in\mbb{N}\qquad\mbox{\it for}\;\;1\leq i\leq n-1,\;0\leq j\leq n-1-i.\eqno(3.83)$$

{\bf Corollary 3.4}. {\it The Verma module} $M_\lmd$ {\it has at most} $n!$ {\it singular vectors up to scalar multiples}. {\it Any singular vector is of the form} $\tau^{-1}(\theta_{\vec i})$ {\it provided} $\theta_{\vec i}$ {\it is a polynomial}. {\it It has exactly} $n!$ {\it singular vectors if and only if} $\lmd$ {\it is a dominant integral weight, that is,} $\lmd_i\in\mbb{N}$ {\it for} $i=1,2....,n-1$. {\it In this case,} $\{\tau^{-1}(\theta_{\vec i})\mid \vec i=(i_1,...,i_{n-1})\in\mbb{N}^{\:n-1},\;i_p\leq p\}$ {\it are all the singular vectors up to scalar multiples}.
\psp

 Finally in this section, we want to list all the polynomial solutions $\theta_{\vec i}$ that involve all $\eta_i$ when $n=3,\;4$. We will not list two solutions that are symmetric under changing indices $i\rta n-i$.

\psp

{\bf Example 3.1}. $n=3$. Polynomial solutions:
$$\eta_1^{\lmd_1+\lmd_2}\eta_2^{\lmd_2+1}(1)=x_{2,1}^{\lmd_1+\lmd_2}x_{3,2}^{\lmd_2+1}\qquad\mbox{with}\;\;\lmd_2\in\mbb{N},\;\;-\lmd_2\leq \lmd_1\in\mbb{Z};\eqno(3.84)$$
\begin{eqnarray*}\hspace{1cm}& &\eta_1^{\lmd_2+1}\eta_2^{\lmd_1+\lmd_2+2}\eta_1^{\lmd_1+1}(1)\\ &=&\sum_{p=1}^{\lmd_1+\lmd_2+2}\frac{\la \lmd_1+1\ra_p\la\lmd_1+\lmd_2+2\ra_p}{p!}x_{2,1}^{\lmd_1+\lmd_2+2-p}x_{3,1}^px_{3,2}^{\lmd_1+\lmd_2+2-p}\hspace{3.1cm}(3.85)\end{eqnarray*}
with $\lmd_1+\lmd_2+1\in\mbb{N}$. 
\pse

{\bf Example 3.2}. $n=4$. Fundamental polynomial solutions are:
\pse

(1) $\eta_3^{\lmd_1+\lmd_2+\lmd_3+3}\eta_2^{\lmd_1+\lmd_2+2}\eta_1^{\lmd_1+1}(1)$ with $\lmd_1\in\mbb{N},\;-\lmd_1-1\leq \lmd_2\in\mbb{Z}$ and $-\lmd_1-\lmd_2-2\leq \lmd_3\in\mbb{Z}$.

(2) $\eta_2^{\lmd_1+\lmd_2+\lmd_3+3}\eta_1^{\lmd_1+1}\eta_3^{\lmd_3+1}(1)$ with $\lmd_1,\lmd_3\in\mbb{N}$ and $-\lmd_1-\lmd_3-2\leq \lmd_2\in\mbb{Z}$.
\pse

(3) $\eta_1^{\lmd_1+\lmd_2+2}\eta_3^{\lmd_2+\lmd_3+2}\eta_2^{\lmd_2+1}(1)$ with $\lmd_2\in\mbb{N},\;-\lmd_2-1\leq \lmd_1\in\mbb{Z}$ and $-\lmd_2-1\leq \lmd_3\in\mbb{Z}$.
\pse

(4) $\eta_1^{\lmd_2+1}\eta_3^{\lmd_1+\lmd_2+\lmd_3+3}\eta_2^{\lmd_1+\lmd_2+2}\eta_1^{\lmd_1+1}(1)$ with $\lmd_1+\lmd_2+1\in\mbb{N}$ and $-\lmd_1-\lmd_2-2\leq \lmd_3\in\mbb{Z}$.

(5) $\eta_2^{\lmd_3+1}\eta_3^{\lmd_1+\lmd_2+\lmd_3+3}\eta_2^{\lmd_1+\lmd_2+2}\eta_1^{\lmd_1+1}(1)$ with $\lmd_1\in\mbb{N}$ and $\lmd_1+\lmd_2+\lmd_3+2\in\mbb{N}$.
\pse

(6) $\eta_1^{\lmd_2+\lmd_3+2}\eta_2^{\lmd_1+\lmd_2+\lmd_3+3}\eta_1^{\lmd_1+1}\eta_3^{\lmd_3+1}(1)$ with $\lmd_3\in\mbb{N}$ and $\lmd_1+\lmd_2+\lmd_3+2\in\mbb{N}$.
\pse

(7)  $\eta_2^{\lmd_1+\lmd_2+\lmd_3+3}\eta_1^{\lmd_1+\lmd_2+2}\eta_3^{\lmd_2+\lmd_3+2}\eta_2^{\lmd_2+1}(1)$ with $\lmd_1+\lmd_2+1,\lmd_2+\lmd_3+1\in\mbb{N}$.
\pse

(8)  $\eta_1^{\lmd_1+\lmd_2+\lmd_3+3}\eta_2^{\lmd_3+1}\eta_3^{\lmd_2+\lmd_3+2}\eta_2^{\lmd_2+1}(1)$ with $\lmd_2+\lmd_3+1\in\mbb{N}$ and $-\lmd_2-\lmd_3-2\leq \lmd_1\in\mbb{Z}$.

(9) $\eta_2^{\lmd_2+\lmd_3+2}\eta_1^{\lmd_2+1}\eta_3^{\lmd_1+\lmd_2+\lmd_3+3}\eta_2^{\lmd_1+\lmd_2+2}\eta_1^{\lmd_1+1}(1)$ with $\lmd_1,\lmd_2,\lmd_3\in\mbb{Z}$ such that $-\lmd_2-1\leq \lmd_1, \lmd_3$ and $-\lmd_2-2\leq \lmd_1+\lmd_2$.
\pse

(10)  $\eta_1^{\lmd_2+\lmd_3+2}\eta_2^{\lmd_3+1}\eta_3^{\lmd_1+\lmd_2+\lmd_3+3}\eta_2^{\lmd_1+\lmd_2+2}\eta_1^{\lmd_1+1}(1)$ with  $\lmd_1+\lmd_2+\lmd_3+2\in\mbb{N}$.
\pse

(11)  $\eta_1^{\lmd_3+1}\eta_2^{\lmd_1+\lmd_2+\lmd_3+3}\eta_1^{\lmd_1+\lmd_2+2}\eta_3^{\lmd_2+\lmd_3+2}\eta_2^{\lmd_2+1}(1)$ with $\lmd_1,\lmd_2,\lmd_3\in\mbb{Z}$ such that $0\leq \lmd_2;\;-\lmd_2-1\leq\lmd_3;-\lmd_2-\lmd_3-2\leq \lmd_1$ or $\lmd_2<-1;\;-\lmd_2-1\leq\lmd_1,\lmd_3$.
\pse

(12)  $\eta_1^{\lmd_3+1}\eta_2^{\lmd_2+\lmd_3+2}\eta_1^{\lmd_2+1}\eta_3^{\lmd_1+\lmd_2+\lmd_3+3}\eta_2^{\lmd_1+\lmd_2+2}\eta_1^{\lmd_1+1}(1)$ with $\lmd_2+1, \lmd_1+\lmd_2+\lmd_3+2\in\mbb{N}$.

\section{Examples of Application}

In this section, we give three examples of applying the results in last section to  the structure of the Verma module $M_\lmd$ for certain special weights $\lmd$ that are not dominant integeral.

The following exchanging relation can be  used to  obtain different presentations of the solutions of the system (2.27) and the singular vectors in the Verma module $M_\lmd$.
\psp

{\bf Proposition 4.1}. {\it For any $\mu_1,\mu_2\in\mbb{C}$ and $1\leq  i<n-1$,  we have}
$$\eta_i^{\mu_1}\eta_{i+1}^{\mu_1+\mu_2}\eta_i^{\mu_2}= \eta_{i+1}^{\mu_2}\eta_i^{\mu_1+\mu_2}\eta_{i+1}^{\mu_1}.\eqno(4.1)$$

{\it Proof}.  Note that for $\mu\in\mbb{C}$, we have
$$[\sum_{p=1}^ix_{i+2,p}\ptl_{i+1,p},x_{i+1,i}^\mu]=\mu x_{i+1,i}^{\mu-1}x_{i+2,i}\eqno(4.2)$$
by (2.40). Moreover,
$$[\sum_{p=1}^ix_{i+2,p}\ptl_{i+1,p},\sum_{j=1}^{i-1}x_{i+1,j}\ptl_{i,j}]=\sum_{j=1}^{i-1}x_{i+2,j}\ptl_{i,j}.\eqno(4.3)$$
Hence
\begin{eqnarray*}& & \eta_i^{\mu_1}\eta_{i+1}^{\mu_1+\mu_2}\\ &=&\sum_{p,q=0}^{\infty}\frac{\la\mu_1\ra_p\la\mu_1+\mu_2\ra_q}{p!q!}x_{i+2,i+1}^{\mu_1+\mu_2-q}x_{i+1,i}^{\mu_1-p}(\sum_{j_1=1}^{i-1}x_{i+1,j_1}\ptl_{i,j_1})^p(\sum_{j_2=1}^ix_{i+2,j_2}\ptl_{i+1,j_2})^q\\ &=& \sum_{p,q,r,s=0}^{\infty}\frac{\la\mu_1\ra_{p+r}\la\mu_1+\mu_2\ra_q\la p\ra_s\la q\ra_{r+s}}{r!s!p!q!}x_{i+2,i+1}^{\mu_1+\mu_2-q}(\sum_{j_2=1}^ix_{i+2,j_2}\ptl_{i+1,j_2})^{q-r-s} x_{i+1,i}^{\mu_1-p-r}\\ & &\times(\sum_{j_1=1}^{i-1}x_{i+1,j_1}\ptl_{i,j_1})^{p-s}
x_{i+2,i}^r(\sum_{j_3=1}^{i-1}x_{i+2,j_3}\ptl_{i,j_3})^s\\ &=&
\sum_{p,q,r,s=0}^{\infty}\frac{\la\mu_1\ra_{p+r}\la\mu_1+\mu_2\ra_q}{r!s!(p-s)!(q-r-s)!}x_{i+2,i+1}^{\mu_1+\mu_2-q}(\sum_{j_2=1}^ix_{i+2,j_2}\ptl_{i+1,j_2})^{q-r-s} x_{i+1,i}^{\mu_1-p-r}\\& &\times (\sum_{j_1=1}^{i-1}x_{i+1,j_1}\ptl_{i,j_1})^{p-s}  x_{i+2,i}^r(\sum_{j_3=1}^{i-1}x_{i+2,j_3}\ptl_{i,j_3})^s\\&=&\sum_{q,k,s=0}^{\infty}\sum_{p=0}^{\infty}\frac{\la\mu_1\ra_k\la\mu_1-k\ra_{ p-s}\la\mu_1+\mu_2\ra_q}{(k-s)!s!(p-s)!(q-k)!}x_{i+2,i+1}^{\mu_1+\mu_2-q}(\sum_{j_2=1}^ix_{i+2,j_2}\ptl_{i+1,j_2})^{q-k} x_{i+1,i}^{\mu_1-k-(p-s)}\\ & &\times(\sum_{j_1=1}^{i-1}x_{i+1,j_1}\ptl_{i,j_1})^{p-s}
x_{i+2,i}^{k-s}(\sum_{j_3=1}^{i-1}x_{i+2,j_3}\ptl_{i,j_3})^s\\ 
&=&\sum_{q,k,s=0}^{\infty}\frac{\la\mu_1\ra_k\la\mu_1+\mu_2\ra_q}{(k-s)!s!(q-k)!}x_{i+2,i+1}^{\mu_1+\mu_2-q}(\sum_{j_2=1}^ix_{i+2,j_2}\ptl_{i+1,j_2})^{q-k} \eta_i^{\mu_1-k}\\ & &\times
x_{i+2,i}^{k-s}(\sum_{j_3=1}^{i-1}x_{i+2,j_3}\ptl_{i,j_3})^s\\&=&\sum_{q,k=0}^{\infty}\sum_{s=0}^{\infty}\frac{\la\mu_1\ra_k\la\mu_1+\mu_2\ra_q}{(k-s)!s!(q-k)!}x_{i+2,i+1}^{\mu_1+\mu_2-q}(\sum_{j_2=1}^ix_{i+2,j_2}\ptl_{i+1,j_2})^{q-k} \eta_i^{\mu_1-k}\\ & &\times
x_{i+2,i}^{k-s}(\sum_{j_3=1}^{i-1}x_{i+2,j_3}\ptl_{i,j_3})^s\\&=&\sum_{q,k=0}^{\infty}\frac{\la\mu_1\ra_k\la\mu_1+\mu_2\ra_q}{k!(q-k)!}x_{i+2,i+1}^{\mu_1+\mu_2-q}(\sum_{j_2=1}^ix_{i+2,j_2}\ptl_{i+1,j_2})^{q-k} \eta_i^{\mu_1-k}(x_{i+2,i}+\sum_{j_3=1}^{i-1}x_{i+2,j_3}\ptl_{i,j_3})^k \hspace{4cm}\end{eqnarray*}
\begin{eqnarray*}   &=&\sum_{k=0}^{\infty}\sum_{q=0}^{\infty}\frac{\la\mu_1\ra_k\la\mu_1+\mu_2\ra_k\la\mu_1+\mu_2-k\ra_{q-k}}{k!(q-k)!}x_{i+2,i+1}^{\mu_1+\mu_2-k-(q-k)}(\sum_{j_2=1}^ix_{i+2,j_2}\ptl_{i+1,j_2})^{q-k} \\ & &\times \eta_i^{\mu_1-k}(x_{i+2,i}+\sum_{j_3=1}^{i-1}x_{i+2,j_3}\ptl_{i,j_3})^k
\\ &=&\sum_{k=0}^{\infty}\frac{\la\mu_1\ra_k\la\mu_1+\mu_2\ra_k}{k!}\eta_{i+1}^{\mu_1+\mu_2-k}\eta_i^{\mu_1-k}(x_{i+2,i}+\sum_{j_3=1}^{i-1}x_{i+2,j_3}\ptl_{i,j_3})^k\hspace{3.7cm}(4.4)\end{eqnarray*}
by (2.36), (4.2) and (4.3). Similarly, we have
$$\eta_i^{\mu_1+\mu_2}\eta_{i+1}^{\mu_1}=\sum_{k=0}^{\infty}\frac{\la\mu_1\ra_k\la\mu_1+\mu_2\ra_k}{k!}\eta_{i+1}^{\mu_1-k}\eta_i^{\mu_1+\mu_2-k}(x_{i+2,i}+\sum_{j_3=1}^{i-1}x_{i+2,j_3}\ptl_{i,j_3})^k.\eqno(4.5)$$
Thus 
\begin{eqnarray*}\eta_i^{\mu_1}\eta_{i+1}^{\mu_1+\mu_2}\eta_i^{\mu_2}&=&
\sum_{k=0}^{\infty}\frac{\la\mu_1\ra_k\la\mu_1+\mu_2\ra_k}{k!}\eta_{i+1}^{\mu_1+\mu_2-k}\eta_i^{\mu_1+\mu_2-k}(x_{i+2,i}+\sum_{j_3=1}^{i-1}x_{i+2,j_3}\ptl_{i,j_3})^k\\ &=& \eta_{i+1}^{\mu_2}\eta_i^{\mu_1+\mu_2}\eta_{i+1}^{\mu_1}.\qquad\Box\hspace{7.9cm}(4.6)\end{eqnarray*}
\pse

{\bf Example 4.1}. Suppose that
$$\lmd_i\in\mbb{N}-1\qquad\for\;\;i=1,2,...,n-1,\eqno(4.7)$$
that is $\lmd_i\in\mbb{N}$ or $\lmd_i=-1$. Set
$$I_0=\{i\in\{1,2,...,n\}\mid \lmd_i\in\mbb{N}\},\;\;I_1=\{1,2,...,n\}\setminus
I_0.\eqno(4.8)$$
Denote by  ${\cal G}_-$  the Lie subalgebra of $sl(n)$ spanned by (2.7) and by $U({\cal G}_-)$ its  universal enveloping algebra. Since 
$$\eta_i^{\lmd_i+1}=\eta^0_i=1\qquad\mbox{if}\;\;\lmd_i=-1,\eqno(4.9)$$
The maximal proper submodule of the Verma module $M_\lmd$  is
$${\cal N}_\lmd=\sum_{i\in I_0}U({\cal G}_-)E_{i+1,i}^{\lmd_i+1}v_\lmd\eqno(4.10)$$
by Corollary 3.4 and induction based on the fact that all $\theta_{\vec i}$ are polynomials.  The character of the irreducible quotient module $V_\lmd=M_\lmd/{\cal N}_\lmd$
can be obtained by using the Weyl character formula for the subalgebra generated by $\{E_{i,i+1}, E_{i+1,i}\mid i\in I_0\}$.
\pse

{\bf Example 4.2}. Suppose that $\lmd_i=-N$ is a negative integer with $N>1$ for some fixed $i\in\{1,2,...,n-1\}$ and
$$\lmd_j\in\mbb{N}-1\qquad\for\;\;i.i-1,i+1\neq j\in\{1,2,...,n-1\}\eqno(4.11)$$
and
$$N-2\leq \lmd_{i-1}\;\;\mbox{if}\;\;i-1>0\;\;\mbox{and}\;\;N-2\leq \lmd_{i+1}\;\;\mbox{if}\;\;i<n-1.\eqno(4.12)$$
Let
$$u=\eta_{r+p}^{p+1+\sum_{j=0}^p\lmd_{r+j}}\eta_{r+p-1}^{p+\sum_{j=0}^{p-1}\lmd_{r+j}}\cdots \eta_r^{\lmd_r+1}(1)\eqno(4.13)$$
such that $r\leq i\leq r+p$. Then
$$\zeta_s(u)=\mu_su\;\;\mbox{and}\;\;\mu_s\in\mbb{N}-1\qquad \for\;\;1\leq s<r+p\eqno(4.14)$$
by (2.30). Assume that 
$$\lmd_{i-1}=N-2\;\;\mbox{if}\;\;i-1>0\;\;\mbox{and}\;\;\lmd_{i+1}=N-2\;\;\mbox{if}\;\;i<n-1.\eqno(4.15)$$
This implies
$$\lmd_{i-1}+\lmd_i+2=\lmd_{i+1}+\lmd_i+2=0.\eqno(4.16)$$
Since $\theta_{\vec i}$ is of form $\eta[\Im]$ (cf. (2.50)-(2.54)), Corollary 3.4 implies that the maximal proper submodule ${\cal N}_\lmd$ of the Verma module $M_\lmd$ is the same as (4.10).
\pse

{\bf Example 4.3}.  Let $N\geq 1$ be an integer. Assume that  $\lmd_i=-N$  for some fixed $i\in\{1,2,...,n-1\}$ and
$$\lmd_j=0\qquad\mbox{if}\;\;j\neq i.\eqno(4.17)$$
Denote by  ${\cal G}_-$  the Lie subalgebra of $sl(n)$ spanned by (2.7) and by $U({\cal G}_-)$ its  universal enveloping algebra. Set
$${\cal N}=\sum_{i\neq j}U({\cal G}_-)E_{j,j+1}v_\lmd.\eqno(4.18)$$
We define the order of the powers in the product formula (3.75) of $\theta_{\vec i}$ starting from ``$(1)$.'' For instance, the first power is $\eta_{i_{n-1}}^{\lmd_{i_{n-1}}+1}$, the $(n-i_{n-1}-1)th$ power is $\eta_{n-2}^{\lmd_{n-1,i_{n-1}}-\lmd_{n-1}-1}$ and the  $(n-i_{n-1})th$ power is $\eta_{n-1}^{\lmd_{n-1,i_{n-1}}}$.
According to (3.62), (3.65) and (3.69)-(3.75), if a negative power of some
 $\eta_j$ with $j\neq i$ appears in $\theta_{\vec i}$, then prior first  power of $\eta_{i-1}$ or  $\eta_{i+1}$ must be negative.  For instance,
$$\eta_{i-1}^{2-N}\eta_{i+1}^{2-n}\eta_i^{1-N}(1)\eqno(4.19)$$
is such a $\theta_{\vec i}$. This implies that the positive powers of $\eta_j$ prior to the first negative power of $\eta_i$ in $\theta_{\vec i}$ do not affect whether $\theta_{\vec i}$ is a polynomial. 

Given $p\geq N$, we have 
$$\eta_{i+p}^{p+1-N}\cdots \eta_{i+N}\eta_{i+N-2}^{-1}\cdots \eta_{i+1}^{2-N}\eta_i^{1-N}(1)=\eta_{i+N-2}^{-1}\cdots \eta_{i+1}^{2-N}\eta_i^{1-N}\eta_{i+p}^{p+1-N}\cdots \eta_{i+N}(1).\eqno(4.20)$$
For  $\mu\in\mbb{C}$, postive integers  $p,j$ such that $j<n-p$, and $p>q\in\mbb{N}$, we have
\begin{eqnarray*}\qquad & &\eta_{j+q}(\eta_{j+p}^{\mu+p}\eta_{j+p-1}^{\mu+p-1}\cdots
\eta_j^{\mu})\\ &=&\eta_{j+p}^{\mu+p}\cdots \eta_{j+q+2}^{\mu+q+2}(\eta_{j+q}\eta_{j+q+1}^{\mu+q+1}\eta_{j+q}^{\mu+q})\eta_{j+q-1}^{\mu+q-1}\cdots
\eta_j^{\mu}\\  &=&\eta_{j+p}^{\mu+p}\cdots \eta_{j+q+2}^{\mu+q+2}(\eta_{j+q+1}^{\mu+q}\eta_{j+q}^{\mu+q+1}\eta_{j+q+1})\eta_{j+q-1}^{\mu+q-1}\cdots\eta_j^{\mu}\\ &=& \eta_{j+p}^{\mu+p}\cdots \eta_{j+q+2}^{\mu+q+2}\eta_{j+q+1}^{\mu+q}\eta_{j+q}^{\mu+q+1}\eta_{j+q-1}^{\mu+q-1}\cdots\eta_j^{\mu}\eta_{j+q+1}\hspace{4.8cm}(4.21)\end{eqnarray*}
by Lemma 4.1. Suppose that $\theta_{\vec i}$ is a polynomial and  
$$\tau(\theta_{\vec i})\not \in {\cal N}+\mbb{C}v_\lmd\eqno(4.22)$$
(cf. (2.22)). By (4.20). (4.21) and induction,  $\theta_{\vec i}$ must be of the form
\begin{eqnarray*}\qquad\theta_{\vec i}&=&\eta_{i+1}^{N-2}\cdots \eta_{i+4-N}\eta_{i+2-N}^{-1}\eta_{i+2}^{N-3}\cdots \eta_{i+6-N}\eta_{i+4-N}^{-1}\eta_{i+3-N}^{-2}\\ & &\cdots
\eta_{i+p-1}^{N-p}\cdots \eta_{i+2p-N}\eta_{i+2p-2-N}^{-1}\cdots\eta_{i+p-N}^{1-p}\cdots\\ &&\eta_{i+N-2} \eta_{i+N-4}^{-1}\cdots \eta_{i-1}^{2-N}\eta_{i+N-2}^{-1}\cdots \eta_{i+1}^{2-N} \eta_i^{1-N}(1)\\ &=&\eta_{i+2-N}^{-1}\eta_{i+3-N}^{-2}\cdots \eta_{i+p-N}^{1-p}\cdots \eta_{i+1}^{2-N} \eta_i^{1-N}(1),\hspace{5.4cm}(4.23)\end{eqnarray*}
which is absurd. Thus the maximal proper submodule ${\cal N}_\lmd$ of the Verma module $M_\lmd$ is ${\cal N}$ in (4.18).
\pse

In a subsequent work, we will show that Examples 4.1 and 4.2 give rise to generalizations of the wedge representations of the Lie algebra $W_{1+\infty}$ (cf. {KP], [KR1]). Moreover, Example 4.3 will imply that the vacuum representations of  $W_{1+\infty}$ with negative integral levels are indeed irreducible.

\hspace{1cm}

\noindent{\Large \bf References}

\hspace{0.5cm}

\begin{description}

\item[{[BGG]}] I. N. Bernstein, I. M. Gel'fand and S. I. Gel'fand,
Structure of representations generated by vectors of highest weight, (Russian) 
{\it Funktsional. Anal. i Prilozhen.} {\bf 8} (1971), no. 1, 1-9.

\item[{[DGK]}] V. V. Deodhar, O. Gabber and V. G. Kac, Structure of some categories of representations of infinite-dimensional Lie algebras, {\it Adv. in Math.} {\bf 45} (1982), no. 1, 92-116.

\item[{[DL]}] V. V. Deodhar and J. Lepowsky, On multiplicity in the Jordan-H\"{o}lder series of Verma modules, {\it J. Algebra} {\bf 49} (1977), no. 2, 512-524.

\item[{[H]}] J. E. Humphreys, {\it Introduction to Lie Algebras and Representation Theory}, Springer-Verlag New York Inc., 1972.

\item[{[J1]}] J. C. Jantzen, Zur charakterformel gewisser darstellungen halbeinfacher grunppen und Lie-algebrun, {\it Math. Z.} {\bf 140} (1974), 127-149.

\item[{[J2]}] J. C. Jantzen, Kontravariante formen auf induzierten Darstellungen habeinfacher Lie-algebren, {\it Math. Ann.} {\bf 226} (1977), no. 1, 53-65.

\item[{[J3]}] J. C. Jantzen, Moduln mit einem h\"{o}chsten gewicht, {\it Lecture Note in Math.} {\bf 750}, Springer, Berlin, 1979.

\item[{[K1]}] V. G. Kac, Highest weight representations of infinite-dimensional Lie algebras, {\it Proc. Intern. Congr. Math. (Helsinki, 1978)}, pp. 299-304, Acad. Sci. Fennica, Helsinki, 1980.

\item[{[K2]}] V. G. Kac,  Some problems on infinite-dimensional Lie algebras and their representations, {\it Lie algebras and related topics (New Brunswick, N.J., 1981)}, pp. 117-126. {\it Lecture Note in Math.} {\bf 933}, Springer, Berlin-New York, 1982.

\item[{[K3]}] V. G. Kac, {\it Infinite-Dimensional Lie algebras}, 3rd Edition, Cambridge University Press, 1990.

\item[{[KK]}] V. G. Kac and D. A. Kazhdan, Structure of representations with highest weight of infinite-dimensional Lie algebras, {\it Adv. in Math.} {\bf 34}
(1979), no. 1, 97-108.

\item[{[KP]}] V. G. Kac and D. H. Peterson, Spin and wedge representations of infinite-dimensional Lie algebras and groups, {\it Proc. Natl. Acad. Sci. USA} {\bf 78} (1981), 3308-3312.

\item[{[KR1]}] V. G. Kac and A. Radul, Quasifinite highest weight modules and the Lie algebra of differential operators on the circle, {\it Commun. Math. Phys.} {\bf 157} (1993), 429-457.

\item[{[KR2]}] V. G. Kac and A. Radul, Representation theory of the vertex algebra $W_{1+\infty}$, {\it Trans. Groups} {\bf 1} (1996), 41-70.

\item[{[L1]}] J. Lepowsky, Canonical vectors in induced modules, {\it Trans. Amer. Math. Soc.} {\bf 208} (1975), 219-272.

\item[{[L2]}] J. Lepowsky, Existence of canonical vectors in induced modules,
{\it Ann. of Math. (2)} {\bf 102} (1975), no. 1, 17-40.

\item[{[L3]}] J. Lepowsky, On the uniqueness of canonical vectors, {\it Proc. Amer. Math. Soc.} {\bf 57} (1976), no. 2, 217-220.

\item[{[L4]}] J. Lepowsky, Generalized Verma modules, the Cartan-Helgason theorem, and the Harish-Chandra homomorphism, {\it J. Algebra} {\bf 49} (1977), no. 2, 470-495.

\item[{[MFF]}] F. G. Malikov,  B. L. Feigin and D. B. Fuchs, Singular vectors in Verma modules over Kac-Moody algebras, (Russian) {\it  Funktsional. Anal. i Prilozhen.} {\bf 20} (1986), no. 2, 25-37.

\item[{[RW1]}]  A. Rocha-Caridi and  N. R. Walllach, Projective modules over graded Lie algebras, I, {\it Math. Z.} {\bf 108} (1982), no. 2, 151-177.

\item[{[RW2]}]  A. Rocha-Caridi and  N. R. Walllach, Highest weight modules over graded Lie algebras, resolutions, filtrations and character formulas,  {\it Trans. Amer. Math. Soc.} {\bf 277} (1983), no. 1, 133-162.

\item[{[S]}] N. N. Sapovolov, A certain bilinear form on the universal enveloping algebra of a complex semisimple Lie algebras, (Russian) {\it  Funktsional. Anal. i Prilozhen.} {\bf 4} (1972), no. 4, 65-70.

\item[{[V1]}] D,-N. Verma, Structure of certain induced representations of cpmlex semisimple Lie algebras, {\it thesis, Yale University,} 1966.

\item[{[V2]}] D,-N. Verma, Structure of certain induced representations of cpmlex semisimple Lie algebras, {\it Bull. Amer. Math. Soc.} {\bf 74}(1968), 160-166.

\item[{[X]}] X. Xu, Differential invariants of classical groups, {\it Duke Math. J.} {\bf 94} (1998), 543-572.

\end{description}

\end{document}